\documentclass[12pt]{article}
\usepackage{amsfonts}

\makeatletter
\newfont{\footsc}{cmcsc10 at 8truept}
\newfont{\footbf}{cmbx10 at 8truept}
\newfont{\footrm}{cmr10 at 10truept}
\renewcommand{\ps@plain}{%
\renewcommand{\@oddfoot}{\footsc  \hfil\footrm\thepage}}
\makeatother 

 \newtheorem{thm}{Theorem}[subsection]
 \newtheorem{cor}[thm]{Corollary}
 \newtheorem{lem}[thm]{Lemma}

\title{Graphs with no $2\delta + 1$ cycle}

\author{Galen E. Turner III\\
\small Mathematics and Statistics Program\\[-0.8ex]
\small Louisiana Tech University, Ruston, Louisiana\\[-0.8ex]
\small \texttt{gturner@coes.LaTech.edu}}

\date{\small Submitted: 9 July 2003\\
\small MR Subject Classifications: 05C35, 05C38, 05C75}

\begin{document}
\maketitle


\begin{abstract}
    Dirac proved that any graph with minimum vertex degree
    $\delta$
    contains either a cycle of length at least $2 \delta$ or a
    Hamilton cycle.  Motivated by this result, we
    characterize those graphs having no cycle longer
    than $2 \delta$.

\end{abstract}

\maketitle

\section{Introduction}

Dirac \cite{Dirac1952} proved that any 2-connected graph with
minimum vertex degree
    $\delta$
    contains either a cycle of length at least $2 \delta$ or a
    Hamilton cycle.
    In this paper, we shall
    characterize those graphs with minimum vertex degree $\delta$
    which have no cycle of length greater than $2 \delta$.
      This characterization was also motivated by Ali
    and Staton \cite{Ali1996} in which similar results were given
    for graphs with no path exceeding $2 \delta + 1$.  In
    particular, as a case in their work, they prove that if the number of vertices of a
    non-hamiltonian 2-connected graph is exactly $2 \delta + 1$, then the graph must
    be isomorphic to the join of a graph on $\delta$ vertices and
    a totally disconnected graph on $\delta + 1$ vertices.  Here,
    we will show that this same type of structure is present when the number of
    vertices exceeds $2 \delta + 1$.


\section{Preliminaries}

    Unless specified otherwise, the terminology used here
    will follow Bondy and Murty \cite{Bon-Mur}.  In particular, a
    graph may have loops and parallel edges.  If $v$ is a
    vertex of a graph $G$, then the set of neighbors of $v$
    will be denoted $N(v)$ or $N_G(v)$, and the minimum vertex
    degree of $G$ will be denoted $\delta$
    or $\delta (G)$.
    If $S$ is a subset of the vertex set of $G$, then
    $G[S]$ denotes the induced subgraph of $G$ on $S$. The
    totally disconnected graph on $m$ vertices will be
    denoted $N_m$, and any subgraph of $K_n$ will be denoted
    $H_n$.
    The {\it join} of two graphs $G$ and $H$ is the graph $G
    \vee H$ and is obtained by taking disjoint copies of $G$
    and $H$ and adding edges joining every vertex of $G$ to
    every vertex of $H$.  In particular, if $|V(G)| = n$ and
    $|V(H)| = m$, then the graph $G \vee H - (E(G) \cup E(H))$
    is isomorphic to $K_{n,m}$.
                Finally, if $x_1 x_2 \ldots x_n$ is a
    path $P$ in a graph $G$, then
    $P[x_i, x_j]$, $P[x_i, x_j)$,
    $P(x_i, x_j]$, and $P(x_i, x_j)$,
    will denote the subpaths $x_i
    x_{i+1} \ldots x_j$, $x_i x_{i+1}
    \ldots x_{j-1}$, $x_{i+1} x_{i+2}
    \ldots x_j$, and $x_{i+1} x_{i+2}
    \ldots x_{j-1}$, respectively.

\section{The Main Theorem}

    To establish the main result of this paper, we shall use
    the following theorem of Dirac \cite{Dirac1952}.

\begin{thm} {\label {Dirac-2delta}}
    Let $G$ be a simple $2$-connected graph with minimum
    vertex degree $\delta$ and suppose $|V(G)| \geq 3$.
    Then $G$ contains either a cycle
    of length at least $2 \delta$ or a Hamilton cycle.
\end{thm}

    Motivated by this result, we will now prove the following
    characterization for the graphs having no cycle
    longer than $2 \delta$.

\begin{thm} {\label {no 2delta+1}}
    Let $G$ be a simple $2$-connected graph with minimum
    vertex degree $\delta$.  If $G$ has
    no cycle of length at least $2 \delta + 1$,
    then $G$ is Hamiltonian or $G \cong H_{\delta} \vee
    N_m$ where $m > \delta$.
\end{thm}

\noindent {\it Proof.}
    Let $G$ be a graph with no Hamilton cycle and no
    cycle of length at least $2 \delta + 1$.
    By Theorem \ref{Dirac-2delta}, G
    has a cycle $C$ of length $2 \delta$.  Now, since $G$ is not Hamiltonian, there is a vertex
    $z$ in $V(G) - V(C)$, and by Menger's Theorem \cite{Menger}, there are 2 paths
    $Q_1$ and $Q_2$ from $z$ to $C$ that have only the vertex $x$
    in common and such that each $Q_i$ meets $C$ in exactly
    one vertex $z_i$.

    \begin{lem} {\label {noname1}}
        $Q_1$ and $Q_2$ must not meet consecutive vertices on $C$.
    \end{lem}

\noindent {\it Proof.}
        Suppose $Q_1$ and $Q_2$ meet $C$ at
         $z_1$ and $z_2$ where $z_1 z_2$ is an edge of $C$.
        Then $G$ has
        a cycle $C'$ formed by the subpath $C - \{z_1z_2\}$ together with
        the path $z_1 z z_{2}$.  Therefore, $C'$ contains
        $2 \delta + 1$ vertices of $G$, namely, all
        the vertices in $V(C) \cup \{z\}$; a contradiction.
    $\Box$  
\vskip .2in

    Since there are paths from every vertex of $V(G) - V(C)$ to
    $C$, choose a longest path $P$ subject to the
    condition that one of the endvertices of $P$ lies on $C$ and no other
    vertex of $P$ lies on $C$.  We label the endvertex of $P$ that is not in $V(C)$ as $x$ and the vertex
     in $V(P) \cap V(C)$ is labelled $y$.  Now, since $P$ is a
     longest path of this type, it is clear that every neighbor of
     $x$ lies in $V(P) \cup V(C)$.  Moreover, we have the
     following lemma.

\begin{lem}{\label{lemma-no-neighbor}}
    $x$ has a neighbor in $V(C) - y$.
\end{lem}

\noindent {\it Proof.}
    Suppose not.  Then every neighbor of $x$ lies on $P$.  Since
    $G$ is $2$-connected, it is clear that there is a vertex of
    $P$ adjacent to some vertex on $C$.  Viewing $P$ as a directed
    path from $x$ to $y$, let $z$ be the first vertex on $P$ that
    is adjacent to a vertex, $z'$, of $C$, and let $u_x$ be the first neighbor of $x$ after $z$ on $P$.
      Now, $C$ is partitioned
    into two paths from $z'$ to $y$, and we arbitrarily label them
    $C_1$ and $C_2$.  Since $C$ has size exactly $2 \delta$, one
    of $C_1$ and $C_2$ has at least $\delta$ edges.  Without loss,
    assume that $C_1$ has at least $\delta$ edges.

        Now, let $P'$ be the path containing $z'z$, $P[z,x]$,
        $xu_x$, and $P[u_x,y]$.  Notice that $P'$ contains
        all the vertices of $P$ except those on $P(z,u_x)$.  Since
        there are no neighbors of $x$ on $P(z,u_x)$, it is clear
        that in addition to $x$ itself, $P'$ contains all the neighbors of
        $x$; thus, $P'$ contains at least $\delta +1$ vertices from $P$ as well as the vertex $z'$ on $C$.
        Since $C_1$ contains $\delta + 1$ vertices, we can combine $C_1$ and $P'$ to obtain a cycle $C'$
         of size $|V(C_1)| + |V(P')| - 2 \geq (\delta + 1) + (\delta + 2)- 2 = 2\delta
        +1$; a contradiction.
$\Box$ \vskip .2in

     Having proved Lemma \ref{lemma-no-neighbor}, let the vertices of $C - y$ that are adjacent
     to $x$ be labelled $x_1, x_2, \ldots,
    x_k$ in cyclic order on $C$, where $y$ lies between $x_k$ and $x_1$
    on $C$, and relabel $y$ as $x_{k+1}$.  It is
    clear that the cycle $C$ is partitioned into $k+1$ internally
    disjoint paths $C(x_1), C(x_2), \ldots, C(x_{k+1})$ where
    $C(x_i)$ is the subpath of $C$ from $x_i$ to $x_{i+1}$ that
    contains no other vertex in $\{x_1, x_2, \ldots, x_{k+1}\}$;
    here $C(x_{k+1})$ is the subpath of $C$ from $x_{k+1}$ to
    $x_1$.
    If $k = 1$, then $C$ can be partitioned into two paths from
    $x_1$ to $x_2$.  We arbitrarily label one of these $C(x_1)$
    and the other $C(x_2)$.  Since $G$ has no cycle containing
    $2\delta + 1$ vertices, it is clear by Lemma \ref{noname1}
    that Each $C(x_i)$ contains at least one vertex in its
    interior, and we let the vertex of $C(x_i)$ that is adjacent to
    $x_i$ be labelled $a_i$.

   \begin{lem} \label{noname3}
        There is no path from $a_i$ to $a_j$ avoiding
        $x \cup C$.
    \end{lem}

    \noindent {\it Proof.}
        Suppose $Q$ is a path joining $a_i$ to $a_j$ and avoiding
        $x \cup C$.  It is easy to see that $G$ has a
        cycle $C'$ formed by the edges of $C - \{a_j x_{j},x_i a_i\}$ together with the path
        $x_jxx_i$ and the edges of $Q$.  Thus, $C'$
        contains the vertices of $C$ and $x$;
         a contradiction.
   $\Box$ 
   \vskip .2in

Now, consider $C(x_{k})$ and $C(x_{k+1})$.  These paths have
    the vertex $x_{k+1}$ as an endvertex.  If one of these paths
    has less than $|P|$ internal vertices, then we can delete the path from
    $C$ and replace it with the path $P$ and the edge $xx_k$ or $xx_{1}$ creating a longer cycle.
    Thus, each of $C(x_k)$ and $C(x_{k+1})$ has length at least
    $|P|+1$.  So, the length of $C(x_j)$ is greater than or equal to $|P|+1$ for $j \in
    \{k,k+1\}$, and by Lemma \ref{noname1}, the length of  $C(x_i)$ is at least $2$ for $i \in
    \{1,\ldots,k-1\}$.  This implies that
    $|C| \geq 2(|P|+1) + (k-1)2$ which means that $ 2\delta \geq 2|P| +
    2k$.  Therefore, $ \delta \geq |P| + k$.

   Now, since $d_G(x) \geq |P|+k$ and $x$ is
    adjacent to exactly $k$ vertices in $V(C) - x_{k+1}$, the vertex $x$ must be adjacent to every vertex on $P-x$.
    This implies that $d_G(x) = |P| + k$ which in turn means that $\delta = |P| +
    k$.  Furthermore, this restriction shows that each $C(x_i)$ has length 2 when $i \in \{1,\ldots
    k-1\}$ and that both $C(x_k)$ and $C(x_{k+1})$ have length
    $|P|+1$.

\begin{lem} \label{noname4}
    $|P| = 1$.  Moreover, $\delta = k + 1$ and $N_G(x) = \{x_1, \ldots, x_{k+1}\}$.
\end{lem}

\noindent {\it Proof.}    Suppose $|P| > 1$ and let the vertices
of $P$ from $x$ to $x_{k+1}$ be labelled $u_1,
    u_2, \ldots, u_{|P|}$ in this order where $u_{|P|} = x_{k+1}$.
    Since $xu_{i+1}$ is an edge, we have a path $P_i = P - \{u_{i}
    u_{i+1}\} \cup xu_{i+1}$ for each $u_i$, and this path is a
    longest path having one endvertex on $C$ and the other
    end being $u_i$.  Since
    $d_G(u_i) \geq \delta$, the vertex $u_i$ must be adjacent to at least $k$ vertices
    on $C - x_{k+1}$, and using Lemmas \ref{noname1} and
    \ref{noname3},
    it is easy to see that the neighbors of $u_i$ on $C - x_{k+1}$ must be
    exactly the vertices
     in $\{x_1,x_2, \ldots, x_{k}\}$.  So, $u_i$
    must be adjacent to every member of $P_i$.  But,
        $(C-\{a_1\}) \cup \{x_2 x u_1 x_1\}$ is a $2\delta + 1$
        cycle, since $u_1 \neq x_{k+1}$.
    $\Box$
    \vskip .2in

   Now, by Lemma \ref{noname4}, every path from a vertex in $V(G) -
   V(C)$ to $V(C)$
   must have length $1$; so, since $G$ is connected, every vertex in $V(G) - V(C)$
   has a neighbor in $V(C)$.
   If $x'$ and $x''$ are
        distinct members of $V(G) - V(C)$, then $x' x''$ is not an edge of $G$; otherwise
        if $y$ is a neighbor of $x'$ on $C$ then $yx'x''$ is a path of length 2 having only one common vertex with $C$.
    We conclude that
    $G[V(G) - V(C)]$ is totally disconnected.

    Upon combining
    Lemmas \ref{noname1}, \ref{noname3}, and  \ref{noname4}, we see that $N(x) = \{x_1, \ldots x_{k+1}\}$
    for every $x$ in $V(G) - V(C)$.  Moreover, the set of neighbors for each
    $a_i$ is $\{x_1, \ldots x_{k+1}\}$ since 
    $a_ia_j$ is not an edge of $G$ 
     and $\delta = k
    +1$.  Therefore,
    $G[\{a_1, \ldots, a_{k+1}\} \cup (V(G) - V(C))]$ is totally disconnected with
    $k+1
    + |V(G) - V(C)| = m$ vertices.  Since $V(G) - V(C)$ has at least one vertex,
    $m$ exceeds $\delta$, and the theorem is established.
$\Box$ \vskip .2in

\begin{cor}
    Let $G$ be a $2$-connected graph with at least $2 \delta + 1$
    vertices.  Then, $G$ has no cycle of length at least $2\delta
    + 1$ if and only if $G \cong H_{\delta} \vee N_m$ where $m >
    \delta$.
\end{cor}


\end{document}